\documentclass[11pt, reqno]{amsart}
\usepackage[top=90pt,bottom=90pt,left=90pt,right=90pt]{geometry}
\usepackage[utf8]{inputenc}
\usepackage[T1]{fontenc}
\usepackage[english]{babel}
\usepackage{enumerate}
\usepackage{fullpage}
\usepackage{bbm,amsmath,amssymb,amsfonts,xcolor,mathtools, verbatim,color,epsf}
\usepackage{bbold}
\usepackage{url}
\usepackage[bookmarksopen=false,breaklinks=true,%
      hyperindex=true,pdfstartview=FitH,colorlinks=true,%
      pdfpagelabels=true,colorlinks=true,linkcolor=black,%
      citecolor=black,urlcolor=blue,hypertexnames=false%
      ]
   {hyperref}
\newcommand{\SSigma}{\mathcal{S}\hspace{-0.5mm}\mathit{\Sigma}}
\newcommand{\N}{\mathbb{N}}
\newcommand{\Z}{\mathbb{Z}}
\newcommand{\R}{\mathbb{R}}
\renewcommand{\S}{\mathcal{S}}
\newcommand{\C}{\mathbb{C}}

\newcommand{\Q}{\mathbb{Q}}

\renewcommand{\leq}{\leqslant}
\usepackage{faktor}
\usepackage{MnSymbol}

\usepackage{etoolbox}
\makeatletter
\patchcmd{\@setaddresses}{\indent}{\noindent}{}{}
\patchcmd{\@setaddresses}{\indent}{\noindent}{}{}
\patchcmd{\@setaddresses}{\indent}{\noindent}{}{}
\makeatother

\newtheorem{theorem}{Theorem}[section]
\newtheorem{proposition}[theorem]{Proposition}
\newtheorem{corollary}[theorem]{Corollary}
\newtheorem{lemma}[theorem]{Lemma}

\newtheorem{remark}[theorem]{Remark}

\DeclareMathOperator\tr{Tr}
\DeclareMathOperator\sper{sper}

\DeclareMathOperator{\sgn}{sgn}

\DeclareFontFamily{U}{mathx}{\hyphenchar\font45}
\DeclareFontShape{U}{mathx}{m}{n}{
      <5> <6> <7> <8> <9> <10>
      <10.95> <12> <14.4> <17.28> <20.74> <24.88>
      mathx10
      }{}
\DeclareSymbolFont{mathx}{U}{mathx}{m}{n}
\DeclareFontSubstitution{U}{mathx}{m}{n}
\DeclareMathAccent{\widecheck}{0}{mathx}{"71}

\title{Hilbert's 17th problem for symmetric functions:\\ undecidability and a counterexample}

\author{Sebastian Debus}
\address{Fachbereich Mathematik und Statistik, Universität Konstanz, 78457 Konstanz, Germany}
\email{sebastian.debus@uni-konstanz.de}

\author{Robin Schabert}
\address{Fachbereich Mathematik und Statistik, Universität Konstanz, 78457 Konstanz, Germany}
\email{robin.schabert@uni-konstanz.de}

\date{\scriptsize{\today}}


\begin{document}

\begin{abstract}
We study inequalities between symmetric functions that hold in every number of variables. First, we prove that it is undecidable whether a symmetric function with rational coefficients is nonnegative in all dimensions, answering a question posed by Acevedo, Blekherman, Debus and Riener. As a consequence, the validity of dimension-free moment criteria for entanglement in quantum information is undecidable. 
Second, we present a homogeneous symmetric function of degree $8$ which is positive definite in every dimension but admits no sum of squares multiplier that is uniform in the number of variables.
Consequently,
Artin's solution of Hilbert's 17th problem and the Krivine-Stengle Positivstellensatz have no
any-dimensional analog and any-dimensional Positivstellensätze asked for by Levin and
Chandrasekaran cannot hold in general.
\end{abstract}

\maketitle

\section{Introduction}
A \emph{symmetric function} $f(X_1,X_2,\dots)$ is a formal power series of bounded degree in countably many variables which is invariant under any permutation of the variables. The ring of symmetric functions, denoted by $\Lambda$, is the polynomial ring $\R[p_1,p_2,\dots]$ in the power sums $p_r = \sum_{i=1}^\infty X_i^r$. Setting all but the first $n$ variables to $0$ truncates each $f \in \Lambda$ to a symmetric polynomial $f^{(n)} \in \R[X_1,\dots,X_n]$. 
Thus a symmetric function $f$ describes a sequence of symmetric polynomials $f^{(n)}$ in an increasing number of variables. An inequality $f \geq 0$ is a statement that must hold in every dimension and is called an \emph{any-dimensional} inequality.
Classical examples are the complete homogeneous symmetric functions of even degree which are positive definite in all dimensions \cite{Hunter1977}. 

The systematic theory of any-dimensional inequalities has recently been studied from three distinct directions: the study of limits of the cones of symmetric nonnegative forms and of sums of squares when the number of variables goes to infinity \cite{acevedo2025symmetric,acevedo2025wonderful,debus2026any}, the framework of any-dimensional polynomial optimization building on de Finetti theorems \cite{levin2025any} and dimension-free trace inequalities which are used as entanglement criteria in quantum information \cite{elben2020mixed, mallick2026optimal, yu2021optimal}.

We call a symmetric function $f \in \Lambda$ \emph{any-dimensional psd} if $f^{(n)} \geq 0$ on $\R^n$ for every $n \in \N$, and $f$ is \emph{any-dimensional sos} if $f^{(n)}$ is a sum of squares for every $n \in \N$. We write $\mathcal{SP}$ for the cone of symmetric functions which are any-dimensional psd and $\S\Sigma$ for the cone of any-dimensional sos. 

Concerning computations, two questions are fundamental. First, is there an algorithm which decides whether a given symmetric function with rational coefficients lies in $\mathcal{SP}$? Second, is there a Positivstellensatz for membership in $\mathcal{SP}$, i.e., a certificate of nonnegativity which is uniform in the number of variables?
For a fixed $n$ the answers to these questions are classical in real algebraic geometry. By the Tarski-Seidenberg theorem \cite{Seidenberg, Tarski} nonnegativity of $f^{(n)}$ is decidable. Moreover, Artin's solution of Hilbert's 17th problem \cite{artin1927zerlegung} proves the existence of a sum of squares $0 \neq q \in \R[X_1,\dots,X_n]$ such that $qf^{(n)}$ is a sum of squares if $f^{(n)}$ is nonnegative. 
Any-dimensional analogs were asked for in \cite[Section~7]{acevedo2025wonderful}, in \cite[Section~8]{debus2026any} and in \cite[Section~7]{levin2025any}. The two main results in this paper are negative answers to both questions.

\subsection*{Main results} We prove in Theorem \ref{thm:undecidable} that the problem of deciding whether a symmetric function $f \in \Lambda$ with rational coefficients is any-dimensional psd is undecidable. The problem remains undecidable even for homogeneous symmetric functions in $\Q[p_2,p_4,\dots,p_{46}]$ and under the promise that either $f \geq \frac{1}{4}p_4^{\frac{\deg f}{4}}$ or $f(x)=-\frac{1}{4}p_4(x)^{\frac{\deg f}{4}}$ for some $x \in \bigcup_{n \in \N} (\R^n \setminus \{0\})$. Thus, deciding any-dimensional positive definiteness is also undecidable.
The theorem is proved by a reduction from Hilbert's 10th problem. The number $46$ comes from Jones \cite{jones1982universal} where it is shown that deciding whether an integer polynomial has a root in $\N_0^9$ is undecidable as opposed to a general $k$ in the DPRM theorem \cite{davis1961decision,matiyasevich1972diophantine}.

We present a homogeneous symmetric function $F$ of degree $8$ which is any-dimensional positive definite. But no formal power series $0 \neq Q \in \R[[X_1,X_2,\dots]]$ of bounded degree exists such that $Q^{(n)}$ and $(QF)^{(n)}$ are sums of squares for all $n \in \N$. This shows that there is no any-dimensional analog of Artin's theorem \cite{artin1927zerlegung} (Theorem \ref{thm:main}). Moreover, no any-dimensional analog of the Krivine-Stengle Positivstellensatz \cite{Krivine,Stengle} exists. It also shows that the hypothesis of the any-dimensional Reznick-type Positivstellensatz in \cite{debus2026any} which requires the truncations to be uniformly bounded below by some $\varepsilon> 0$ on the unit sphere cannot be weakened to positive definiteness.

Finally, in Section \ref{sec:quantum} we translate Theorems \ref{thm:undecidable} and \ref{thm:main} into the language of \cite{mallick2026optimal}. 
We obtain that the validity of dimension-free trace inequality criteria for separability, in particular criteria based on moments of the partial transpose, cannot be decided algorithmically (Theorem \ref{thm:quantum1}).

\subsection*{Related work and comparison with power mean inequalities}
Related but distinct work concerns any-dimensional inequalities in \emph{power means} \cite{AcevedoBlekherman,BlekhermanRiener}, or equivalently in normalized traces and moments \cite{Klep26,Klep21}. In this framework a symmetric function $f \in \Lambda$, expressed as $f=g(p_1,\dots,p_d)$ for some unique polynomial $g \in \R[Y_1,\dots,Y_d]$, is not truncated to $f^{(n)}$ but the truncations $p_r^{(n)}$ are replaced by the power means $\frac{1}{n} p_r^{(n)}$. 
Since for $x \in \R^n$ the evaluation $\frac{1}{n}p_r^{(n)}(x)$ is the $r$-th moment of the discrete probability measure $\sum_{i=1}^n \frac{1}{n} \delta_{x_i}$, nonnegativity in all dimensions actually corresponds to nonnegativity of the polynomial $g$ on the Hankel spectrahedron of truncated moment sequences of probability measures on the real line \cite[Proposition 4.1]{AcevedoBlekherman}. Consequently nonnegativity in this framework is decidable by real quantifier elimination. In the language of trace polynomials with normalized trace there exists a Krivine-Stengle Positivstellensatz \cite{Klep21} which was also extended to several variables through moment polynomials \cite{Klep26}. 

In the setting investigated in this paper the probability measures $\sum_{i=1}^n \frac{1}{n} \delta_{x_i}$ are replaced by the counting measures $\sum_{i=1}^n \delta_{x_i}$ with integer weights. 
It follows from Timofte's half-degree principle \cite{Riener12,Timofte} that nonnegativity of $f^{(n)}$ is equivalent to nonnegativity of $f^{(n)}$ at points with few distinct coordinates whose multiplicities are integers summing to $n$. This countable discrete structure is visible in the infinitely many extreme rays of the degree three cone of dimension-free trace inequalities \cite[Proposition 3.5]{mallick2026optimal} and in the characterization of the character space of the preordering of any-dimensional sums of squares in \cite[Corollary 4.9]{debus2026any}. Moreover, in \cite{acevedo2025wonderful} the authors proved in the spirit of Hatami and Norine's undecidability theorem for graph homomorphism densities \cite{hatami2011undecidability} that any-dimensional nonnegativity is undecidable for symmetric functions in several groups of distinct variables. The case of a single group of variables, i.e., the elements of $\mathcal{SP}$, was left open in \cite{acevedo2025wonderful}. While finalizing this paper we learned of the manuscript \cite{Jin2026} which independently settles this case (see Remark \ref{rem:Jin} for a comparison).

Moreover, we note that the undecidability result (Theorem \ref{thm:undecidable}) follows also from a more general statement on undecidability of nonnegativity on certain non-semialgebraic compact sets \cite{Riener}.

\subsection*{Notation}
We write $\N=\{1,2,3,\dots\}$ and $\N_0 = \N \cup \{0\}$. The ring of symmetric functions is $\Lambda=\R[p_1,p_2,\dots]$ with $p_r=\sum_{i=1}^\infty X_i^r$ and $\Lambda_\Q=\Q[p_1,p_2,\dots]$ is the ring of symmetric functions over the rationals. The truncation of $f \in \Lambda$ to $n$ variables is $f^{(n)}(X_1,\dots,X_n)=f(X_1,\dots,X_n,0,0,\dots)$ and is therefore  determined by $p_r^{(n)}=X_1^r+\dots+X_n^r$. We sometimes write $f(X_1,\dots,X_n)$ for $f^{(n)}(X_1,\dots,X_n)$. 
For a formal power series $g \in \R[[X_1,X_2,\dots]]$ we write $g^{(n)}$ for its truncation to the first $n$ variables, i.e., $g^{(n)}=g(X_1,\dots,X_n,0,0,\dots)$. The degree of $g$ is defined as $\deg g = \sup_{n \in \N} \deg g^{(n)}$ and $g$ is of bounded degree if $\deg g < \infty$. 
A symmetric function $f$ is homogeneous of degree $d$ if it is a linear combination of products $p_{r_1}\dots p_{r_s}$ with $r_1+\dots+r_s=d$. Then every truncation $f^{(n)}$ is a homogeneous polynomial of degree $d$ or zero. We write $\Lambda_d$ for the homogeneous part of degree $d$ in $\Lambda$. A polynomial $f$ is even if all of its monomials have only even exponents. We write $\Lambda^{\varepsilon}:=\R[p_2,p_4,\dots]$ for the set of even symmetric functions. All norms used below are Euclidean. Finally, a polynomial is sos if it is a sum of squares of polynomials.

\section{Undecidability}\label{sec:undecidability}
In this section, we prove that deciding nonnegativity of symmetric functions with rational coefficients in all numbers of variables is undecidable. 
\begin{remark}\label{rem:undecidability}
It follows from the theorem of Davis, Putnam, Robinson and Matiyasevich \cite{davis1961decision,matiyasevich1972diophantine} that deciding whether a Diophantine equation has a solution in $\N_0^k$ is undecidable, which is the negative solution of Hilbert's 10th problem. Jones \cite{jones1982universal} proved that deciding whether a Diophantine equation has a solution in $\N_0^9$ is already an undecidable problem which was announced earlier by Matiyasevich \cite{matijasevivc1977some}. 
\end{remark}

In \cite[Theorem~6.2]{acevedo2025wonderful} it was shown that the problem of deciding nonnegativity of symmetric functions in several groups of countably many disjoint variables is undecidable. The proof was inspired by Hatami and Norine \cite{hatami2011undecidability} where undecidability for inequalities in graph homomorphism densities was proven. In \cite{acevedo2025wonderful} the authors used the Vandermonde cell $\Pi_3$, i.e., the Euclidean closure of the set $\bigcup_{n \in \N} (p_2,p_3)(\{ x \in \R_{\geq 0}^n : \sum_{i=1}^n x_i = 1\})$, to relate nonnegativity of polynomials on $\Pi_3^k$ to nonnegativity of integer polynomials in $k$ variables.

Here we encode integer polynomials $P$ by homogeneous symmetric functions $H_P$ and use a rounding procedure to relate the property of $P$ having a root in $\N_0^k$ to any-dimensional positive definiteness of $H_P$.

\begin{theorem}\label{thm:undecidable}
    The following decision problem is undecidable:
    \begin{enumerate}
        \item[] Instance: A symmetric function $f(X_1,X_2,\dots)$ with rational coefficients.
        \item[] Question: Is the truncation of $f(X_1,X_2,\dots)$ to $n$ variables nonnegative on $\Q^n$ for all $n \in \N$?
    \end{enumerate}
    The problem remains undecidable when the instances are restricted to homogeneous even symmetric functions in $\Q[p_2,p_4,\dots,p_{46}]$.
\end{theorem}
Since each truncation $f^{(n)}$ is a polynomial and thus continuous, the question in Theorem 2.2 is equivalent to asking whether $f$ is any-dimensional psd. 

Moreover, the theorem shows that no effective Positivstellensätze for any-dimensional polynomial optimization problems exist in general. This is because the set $\mathcal{SP} \cap \Lambda_\Q$ is not recursively enumerable and so no checkable certificate system can be complete. 

In the language of \cite{acevedo2025wonderful} we obtain that deciding nonnegativity of rational polynomials on the Vandermonde cell 
\[\Pi_{23}:=\operatorname{cl} \bigcup_{n \in \N} (p_2,\dots,p_{23})\left(\left\{x \in \R^n_{\geq 0} ~\middle|~ \sum_{i=1}^n x_i=1\right\}\right)\] is undecidable. This follows from Theorem \ref{thm:undecidable}, since $H_P \in \Q[p_2,p_4,\dots,p_{46}]$ is homogeneous and so nonnegativity is equivalent to nonnegativity on the sphere, and from substituting $x_i$ by $\sqrt{x}_i$.

\subsection{Strategy of the proof of Theorem \ref{thm:undecidable}}

For $P \in \Z[Y_1,\dots,Y_k]$ we compute a symmetric function $H_P \in \Q[p_2,p_4,\dots,p_{4k+10}]$ which is any-dimensional psd if and only if $P$ has no zero in $\N_0^k$. The undecidability then follows from Remark \ref{rem:undecidability}.

Elements in $\N_0^k$ enter as multiplicity vectors. For fixed labels $a_1,\dots,a_k \in \Q$ a vector with coordinates in $\{0,a_1,\dots,a_k\}$ will be associated with the vector $N \in \N_0^k$ counting the multiplicities of each label. Every $N \in \N_0^k$ occurs for some vector length.
We then use Lagrange polynomials $\ell_j \in \Q[T]$ with $\ell_j(0)=0$ and $\ell_j(a_i)=\delta_{ij}$ to form the symmetric functions $S_j = \sum_i \ell_j (X_i)$ and $S:=(S_1,\dots,S_k)$. Then $P(S)^2-\frac{1}{4}$ equals $-\frac{1}{4}$ on every vector such that the associated $N \in \N_0^k$ is a zero of $P$. On all other configurations whose coordinates are labels we have $P(S)^2-\frac{1}{4}\geq \frac{3}{4}.$
We bound the rounding error $||S-N||^2$ in Lemma \ref{lem:rounding} through a penalty term which is also a symmetric function. The penalty term vanishes on $\bigcup_{n \in \N}\R_{\geq 0}^n$ precisely on configurations whose coordinates lie in $\{0\} \cup \mathcal{A}$. This allows us to compute a symmetric function $Q_P$ which satisfies $Q_P \geq \frac{1}{4}$ on $\bigcup_{n \in \N}\R_{\geq 0}^n$ if and only if $P$ has no zero in $\N_0^k$ (Lemma \ref{lem:encoding}). 
Finally, a substitution of variables removes the sign condition and homogenizes $Q_P$ and we obtain $H_P \in \Lambda_\Q$. 

\subsection{Approximating integer vectors with symmetric functions}
Throughout the remainder of this section we fix a positive integer $k$, define 
\[\mathcal{A}:=\left\{\frac{1}{2},1,2,\dots,k+1\right\}, \quad m := k+2 \quad \text{and} \quad a_j := j+1 \quad \text{for }  1 \leq j \leq k.\]
Furthermore, we consider the univariate rational polynomials \[Q(T):=\prod_{a \in \mathcal{A}} (T-a) \quad \text{and} \quad R(T):=T\cdot Q(T).\]
Thus, $\mathcal{A}$ contains $m$ elements and $\{0\} \cup \mathcal{A}$ contains $m+1$ elements. Note that any two distinct elements in $\{0\} \cup \mathcal{A}$ have distance at least $\frac{1}{2}$. 
The numbers $a_1,\dots,a_k$ will be used as labels. Multiplicities of coordinates in vectors from $(\{0\}\cup \mathcal{A})^n$ will be used to encode vectors in $\N_0^k$ and the additional labels $\frac{1}{2}$ and $1$ will allow a homogenization.

Moreover, for $1 \leq j \leq k$ we consider the polynomials \[ \ell_j(T) := \frac{R(T)}{(T-a_j) \cdot R'(a_j)} = \frac{T}{a_j} \prod_{a \in \mathcal{A} \setminus \{a_j\}}  \frac{T-a}{a_j-a} \in \mathbb{Q}[T].\]
Observe that $\deg \ell_j = m$, $\ell_j(0)=0$ and $\ell_j(a)= \delta_{a,a_j}$ for $a \in \mathcal{A}$, where $\delta_{p,q}$ denotes the Kronecker delta function. 
Thus, $\ell_1,\dots,\ell_k$ are Lagrange polynomials.
We also note $R'(a_j) = a_j \prod_{a \in \mathcal{A}\setminus \{a_j\}} (a_j-a) \neq 0$. Finally, for $t \geq 0$ let $\rho(t)$ denote the number in $\{0\} \cup \mathcal{A}$ closest to $t$, where the smaller number is chosen in a tie, and we define $\rho_j(t):=\delta_{\rho(t),a_j}$ for $1 \leq j \leq k$.
\begin{lemma}\label{lem:rounding}
There exists an effectively computable positive integer $C$, depending only on $k$, for which
\begin{align}\label{eq:1}
    |\ell_j(t)-\rho_j(t)| \leq C\cdot t \cdot |Q(t)| \qquad \forall t \geq 0,\, 1 \leq j \leq k.
\end{align}
Moreover, for every $n \in \N$ and $x=(x_1,\dots,x_n) \in \R_{\geq 0}^n$, set
\begin{align}\label{eq: S, N, U, E}
S_j := \sum_{i=1}^n \ell_j(x_i), \quad N_j := \sum_{i=1}^n \rho_j(x_i), \quad U:=\sum_{i=1}^n x_i,\quad  E := \sum_{i=1}^n x_i Q(x_i)^2 \quad \text{and} \quad S:=(S_1,\dots,S_k).
\end{align}
Then $N=(N_1,\dots,N_k)\in \N_0^k$ and 
\begin{align}\label{eq:2}
    ||S-N||^2 \leq k \cdot C^2 \cdot U\cdot E, \qquad ||N||^2 \leq 16\cdot U^2.
\end{align}
\end{lemma}
Observe that $U=p_1(x)$ and $S_j(x),E(x)$ are symmetric in $x$ and can be expressed in the power sum basis independently of $n$ (see Remark \ref{rem: S,U,E}). So, we can control the distance between $N \in \N_0^k$ and the symmetric quantities in \eqref{eq:2} by evaluating symmetric polynomials.
\begin{proof}
Note that for $t \geq 0$ and $b \in \{0\} \cup \mathcal{A}$ with $\rho(t) \neq b$ we have $|t-b|\geq \frac{1}{4}$, since the minimal distance between numbers from $\{0\} \cup \mathcal{A}$ is $\frac{1}{2}$.

Since $a_j$ is a root of $\ell_j(T)-1$ we find that $M_j(T) := \frac{\ell_j(T)-1}{T-a_j} \in \Q[T]$ is a polynomial of degree $m-1$. Let $B_j$ denote the sum of the absolute values of the coefficients of $M_j$ and set $T_0 := 2(k+1)$. 
We claim that any integer
\begin{align}\label{eq: C bound}
    C \geq \max \left\{ 1,\max_{1 \leq j \leq k} \frac{4}{|R'(a_j)|}, \max_{1 \leq j \leq k}4^m B_j T_0^{m-1} \right\}
\end{align}
satisfies \eqref{eq:1}.
Since all occurring numbers in the maximum are rational and can be computed from $k$ in finitely many steps, the integer $C$ is effectively computable.

To verify \eqref{eq:1} we fix $t \geq 0$ and $1 \leq j \leq k$. We proceed by a case distinction on whether $\rho(t)=a_j$.

If $\rho(t)\neq a_j$ then $\rho_j(t)=0$ and $|t-a_j| \geq \frac{1}{4}$. We obtain
\[ |\ell_j(t)-\rho_j(t)|=|\ell_j(t)|=\frac{|R(t)|}{|t-a_j| \cdot |R'(a_j)|} \leq \frac{4}{|R'(a_j)|} \cdot t \cdot |Q(t)| \leq C\cdot t\cdot |Q(t)|\]
which proves the desired estimate.

If $\rho(t) = a_j$ we have $\rho_j(t)=1$ and $\ell_j(t)-\rho_j(t) = \ell_j(t)-1=(t-a_j)M_j(t)$. 
Consider the monic polynomial
\[D_j(T) := \frac{R(T)}{T-a_j} = \prod_{ b \in (\{0\} \cup \mathcal{A}) \setminus \{a_j\}} (T-b)\] of degree $m$. 
Since for every root $b$ of $D_j(T)$ we have $|t-b| \geq \frac{1}{4}$, we also have $|D_j(t)| \geq \frac{1}{4^m}$. 
Next, we show that 
\begin{align}\label{eq:3}
 |M_j(t)| \leq 4^m B_j T_0^{m-1} |D_j(t)|.   
\end{align}
Note that for all $t \geq 0$ we have $|M_j(t)| \leq B_j \max\{1,t\}^{m-1}$ since $M_j$ has degree $m-1$ and $B_j$ is the sum of the absolute values of its coefficients.

If $0 \leq t \leq T_0$ then $|M_j(t)| \leq B_j \cdot \max \{1,t\}^{m-1} \leq B_j T_0^{m-1}$ which implies \eqref{eq:3}. Since $T_0 \geq 1$ by definition, we also get $\max\{1,t\}^{m-1} \leq T_0^{m-1}$.

If $t \geq T_0$ then for every root $b$ of $D_j(T)$ we have $t-b \geq t-(k+1) \geq \frac{t}{2} > 0$, since $b \in [0,k+1]$. Thus, $|D_j(t)| \geq \left(\frac{t}{2}\right)^m$. With $|M_j(t)|\leq B_jt^{m-1}$ we obtain
\[ \frac{|M_j(t)|}{|D_j(t)|} \leq \frac{2^mB_j}{t} \leq \frac{2^m B_j}{T_0} \leq 4^m B_j T_0^{m-1},\]
where the last inequality follows from $T_0 \geq 1$. Then multiplication of $|M_j(t)| \leq 4^m B_j T_0^{m-1} |D_j(t)|$ by $|t-a_j|$ yields
\[ |\ell_j(t)-\rho_j(t)| = |\ell_j(t)-1| \leq 4^m B_j T_0^{m-1}|t-a_j| \cdot |D_j(t)| = 4^m B_j T_0^{m-1} |R(t)| \leq Ct|Q(t)| \]
which completes the proof of \eqref{eq:1}.

Now let $n \in \N$ and $x \in \R_{\geq 0}^n$. By construction we have $N \in \N_0^k$. We use \eqref{eq:1} and apply the Cauchy-Schwarz inequality to the vectors $(\sqrt{x_i})_i, (\sqrt{x_i}\cdot|Q(x_i)|)_i \in \R^n$. Then
\begin{align*}
    |S_j-N_j| \leq \sum_{i=1}^n |\ell_j(x_i)-\rho_j(x_i)| \leq C \sum_{i=1}^n x_i |Q(x_i)| \leq C \left( \sum_{i=1}^n x_i\right)^{\frac{1}{2}} \left( \sum_{i=1}^n x_iQ(x_i)^2\right)^{\frac{1}{2}} = C\sqrt{UE}.
\end{align*}
The first inequality in \eqref{eq:2} follows now by squaring and summing over $1 \leq j \leq k$.
To see the second inequality, note that
if $\rho(x_i) \neq 0$ then $ \rho(x_i) \geq \frac{1}{2}$ and $\rho(x_i)$ is the closest point to $x_i$ in $\{0\} \cup \mathcal{A}$. Thus $|x_i-\rho(x_i)| \leq |x_i-0|=x_i$. Since $x_i = \rho(x_i)+(x_i-\rho(x_i)) \geq \rho(x_i) - |x_i-\rho(x_i)| \geq \frac{1}{2}-x_i$ we obtain $x_i \geq \frac{1}{4}$.
Hence,
\[ \sum_{j=1}^k N_j = \# \left\{i \in [n] ~\middle|~ \rho(x_i) \in \left\{a_1,\dots,a_k\right\}\right\} \leq \# \left\{i \in [n] ~\middle|~ \rho(x_i) \neq 0 \right\} \leq \sum_{\substack{i \in [n] \\ \rho(x_i) \neq 0}} 4x_i \leq 4U \]
which implies $||N||^2 \leq (\sum_{j=1}^k N_j)^2 \leq 16U^2$.
\end{proof}

Note that \eqref{eq: C bound} gives an effectively computable lower bound for $C$.

\begin{remark}\label{rem: S,U,E}
    Since the polynomials $\ell_j(T)$ and $T \cdot Q(T)^2$ depend only on $k$, we can write 
    \[\ell_j(T)=\sum_{r=1}^m c_{j,r}T^r \quad \text{and} \quad T\cdot Q(T)^2 = \sum_{r=1}^{2m+1}q_rT^r,\] 
    where $c_{j,r},q_r \in \Q$. Note that $\ell_j$ has no constant term since $\ell_j(0)=0$. For $x \in \bigcup_{n \in \N}\R^n_{\geq 0}$ and $S_j$, $U$ and $E$ as in Lemma \ref{lem:rounding}, we can therefore write
    \[ S_j = \sum_{r=1}^m c_{j,r} p_r(x), \quad U = p_1(x) \quad \text{and} \quad E=\sum_{r=1}^{2m+1}q_rp_r(x).\]
    So we can and will consider $S_j$, $U$ and $E$ as symmetric functions in $\Q[p_1,\dots,p_{2m+1}]$ whose restrictions to $\bigcup_{n \in \N}\R_{\geq 0}^n$ satisfy \eqref{eq:2} in Lemma \ref{lem:rounding}.
\end{remark}

\subsection{An encoding of integer polynomials}
Define $P(Y) := \sum_{\alpha \in \N_0^k} c_\alpha Y^\alpha \in \Z[Y_1,\dots,Y_k]$, let $C$ be as in Lemma \ref{lem:rounding} and $d :=\max\{1,\deg P\}$, where we use the convention $\deg 0 = -\infty$. Define
\[ b_j := \sum_\alpha |c_\alpha| \alpha_j \in \N_0 \quad \text{and} \quad K_P := 1+kC^2 \sum_{j=1}^kb_j^2 \in \N \quad \text{for $1 \leq j \leq k$}.\]
So, $K_P$ is a positive effectively computable integer depending on the polynomial $P$ and $k$.
With $S=(S_1,\dots,S_k),U,E \in \Q[p_1,\dots,p_{2m+1}]$ as in Remark \ref{rem: S,U,E} we define the symmetric function
\begin{align}\label{eq:4}
    Q_P := P(S_1,\dots,S_k)^2 + K_PUE \left( 1+ \sum_{j=1}^k S_j^2+ 16U^2\right)^{d-1}-\frac{1}{4}.
\end{align}
It follows from Remark \ref{rem: S,U,E} that
\begin{align}\label{eq:degreebound}
    \deg S_j \leq m, \quad \deg U = 1 \quad \text{and} \quad \deg E \leq 2m+1
\end{align}
which implies $\deg Q_P \leq L:=2md+2$.
\begin{lemma}\label{lem:encoding}
    Let $P \in \Z[Y_1,\dots,Y_k]$.
    \begin{enumerate}
        \item[(a)] If $P$ has no zero in $\N_0^k$, then $Q_P(x) \geq \frac{1}{4}$ for every $x \in \bigcup_{n \in \N} \R_{\geq 0}^n$. 
        \item[(b)] If $P(N)=0$ for some $N \in \N_0^k$ then there exist $n \in \N$ and $z \in \mathcal{A}^n$ such that
        \[ p_1(z)=p_2(z) \qquad \text{and} \qquad Q_P(z)=-\frac{1}{4}.\] 
    \end{enumerate}
\end{lemma}
\begin{proof}
    We write $V := 1 + \sum_{j=1}^k S_j^2+16U^2 = 1 + ||S||^2 + 16 U^2 \geq 1$ and so $Q_P=P(S)^2+K_PUEV^{d-1}-\frac{1}{4}$.
    \begin{enumerate}
        \item[(a)] Let $x \in \R_{\geq 0}^n$ for some $n \in \N$ and let $N \in \N_0^k$ be the vector from Lemma \ref{lem:rounding}. 
         We first prove the auxiliary claim
    \begin{align}\label{eq:6}
        |P(S)-P(N)|^2  \leq K_PUEV^{d-1},
    \end{align}
    where the right-hand side is the middle term of $Q_P$.
        For $\theta \in [0,1]$ we set $\xi(\theta):= N+\theta(S-N)$ and write $\xi := \xi(\theta)$.
        Since the Euclidean norm is convex and since $||N||^2\leq 16U^2$ by \eqref{eq:2} in Lemma \ref{lem:rounding}, we obtain
        \[ 1 + ||\xi(\theta)||^2 \leq 1 + \max\{||S||^2,||N||^2\} \leq 1+||S||^2+16U^2=V.\]
        For $1 \leq j \leq k$ every monomial appearing in the partial derivative $\partial_j P = \sum_{\substack{\alpha\in \N_0^k \\ \alpha_j \geq 1}}  c_\alpha \alpha_j Y^{\alpha-e_j}$ has degree at most $d-1$. We have $|\xi_i| \leq \sqrt{1+||\xi||^2}$ and $1 + ||\xi||^2 \leq V$. Thus, every monomial $Y^\beta$ with $|\beta| \leq d-1$ satisfies $|\xi^\beta| \leq (1+||\xi||^2)^\frac{|\beta|}{2} \leq V^{\frac{d-1}{2}}$. Thus the inequality
        \[ | \partial_j P(\xi(\theta))| \leq b_j V^{\frac{d-1}{2}} \qquad \text{for all }0 \leq \theta \leq 1\text{ and } 1 \leq j \leq k,\]
        holds.
        Since $P(S)-P(N)=\int_0^1 \langle \nabla P(\xi(\theta)), S-N\rangle \mathrm{d}\theta,$ we obtain from the Cauchy-Schwarz inequality and \eqref{eq:2}
        \[ |P(S)-P(N)|^2 \leq \left( \sum_{j=1}^k b_j^2\right)V^{d-1}||S-N||^2 \leq \left( \sum_{j=1}^k b_j^2 \right) V^{d-1}kC^2 UE \leq K_PUEV^{d-1}\]
        which is precisely the claimed inequality \eqref{eq:6}. 
        
        Finally, note that $|P(N)| \geq 1$ since $P$ has no zero in $\N_0^k$. Thus,
        \[ 1 \leq |P(N)|^2 \leq 2P(S)^2 +2 |P(S)-P(N)|^2 \leq 2\left(P(S)^2+K_PUEV^{d-1}\right), \]
        and so 
        \[Q_P(x) = P(S)^2+K_PUEV^{d-1}-\frac{1}{4} \geq \frac{1}{2}-\frac{1}{4}= \frac{1}{4}.\]
        \item[(b)] Next, we suppose $P(N)=0$ for some $N \in \N_0^k$. We consider the vector $z$ consisting of $N_j$ copies of $a_j=j+1$ for each $1 \leq j \leq k$, $B:=4\sum_{j=1}^k a_j (a_j-1)N_j \in \N_0$ copies of $\frac{1}{2}$ and one copy of $1$. By construction we have $z \in \mathcal{A}^n$, where $n = \sum_j N_j+B+1 \geq 1.$
        Every coordinate of $z$ lies in $\mathcal{A}$ which implies $E=0$ since $Q(z_i)=0$ for all $1 \leq i \leq n$. Since $\ell_j(a)=\delta_{a,a_j}$ for all $ a \in \mathcal{A}$ and $\frac{1}{2},1 \not\in \{a_1,\dots,a_k\}$ we obtain $S_j=N_j$ for all $j$ and so $P(S)=P(N)=0$ and $Q_P(z)=-\frac{1}{4}$.
        Moreover,
        \[ p_2(z)-p_1(z)=\sum_{i=1}^n z_i(z_i-1)=\sum_{j=1}^k a_j(a_j-1)N_j+B\frac{1}{2}\left(\frac{1}{2}-1\right)+1 \cdot 0 = \sum_{j=1}^k a_j(a_j-1)N_j-\frac{B}{4}=0\]
        and thus $p_1(z)=p_2(z)$ as claimed.
    \end{enumerate}
\end{proof}

\subsection{Proof of Theorem \ref{thm:undecidable}}
For $P \in \Z[Y_1,\dots,Y_k]$ the symmetric function $Q_P$ defined in \eqref{eq:4} (see also Remark \ref{rem: S,U,E}) is not necessarily homogeneous and Lemmas \ref{lem:rounding} and \ref{lem:encoding} relate zeros of $P$ in $\N_0^k$ to evaluations of $Q_P$ on the set $\bigcup_{n \in \N}\R_{\geq 0}^n$ instead of $\bigcup_{n \in \N}\R^n$. We resolve this by substituting each power sum function $p_r$ by $\frac{p_2^rp_{2r}}{p_4^r}$ for $1 \leq r \leq 2m+1$ in $Q_P$ and multiplying by $p_4^L$. The resulting symmetric function $H_P$ is homogeneous and even. 

For $n \in \N$ and $0 \neq x \in \R^n$ we set
\begin{align}\label{eq:a,b,z}
   y=(y_1,\dots,y_n):= \left(\frac{p_2(x)}{p_4(x)}x_1^2,\dots,\frac{p_2(x)}{p_4(x)}x_n^2\right) \in \R_{\geq 0}^n.
\end{align}
Since we assumed $0 \neq x$ we have $p_2(x),p_4(x) > 0$. 
Then for any $r \geq 1$ we have
\[ p_r(y)=\frac{p_2(x)^r}{p_4(x)^r}p_{2r}(x)\]
and $p_1(y)=\frac{p_2(x)^2}{p_4(x)}=\frac{p_2(x)^2}{p_4(x)^2}p_4(x) = p_2(y).$

Recall that $Q_P \in \Q[p_1,\dots,p_{2m+1}]$ for $P \in \Z[Y_1,\dots,Y_k]$. So there exists a unique polynomial $G_P \in \Q[Z_1,\dots,Z_{2m+1}]$ with $Q_P=G_P(p_1,\dots,p_{2m+1})$. 
We define the symmetric function
\begin{align}\label{eq:7}
    H_P := p_4^L G_P \left( \frac{p_2^2}{p_4},\frac{p_2^2 p_4}{p_4^2}, \frac{p_2^3 p_6}{p_4^3}, \dots, \frac{p_2^{2m+1}p_{4m+2}}{p_4^{2m+1}}\right).
\end{align}
Formally, we obtain $H_P$ from $Q_P$ by substituting each $p_r$ by $\frac{p_2^rp_{2r}}{p_4^r}$ for $1 \leq r \leq 2m+1$ and multiplying by $p_4^L$ to clear the denominators. The computations are done in the quotient field of $\Q[p_2,p_4,p_6,\dots]$. Then $H_P$ is indeed a homogeneous even symmetric function (see Proposition \ref{prop:undecidable}). 
Note that the first two substituted entries in the evaluation of $G_P$ in \eqref{eq:7} both equal $\frac{p_2^2}{p_4}$ and therefore coincide. This happens intentionally and reflects the normalization $p_1(z)=p_2(z)$ in Lemma \ref{lem:encoding} (b). 

\begin{proposition}\label{prop:undecidable}
For $P \in \Z[Y_1,\dots,Y_k]$ we have $H_P \in \Q[p_2,p_4,\dots,p_{4m+2}]$ and $H_P$ is computable by finite rational arithmetic from $P$. Moreover, $H_P$ is a homogeneous even symmetric function of degree $4L=8(k+2)d+8$ and satisfies the following properties:
\begin{enumerate}
    \item[(a)] If $P$ has no zero in $\N_0^k$, then for every $n \in \N$ and $0 \neq x \in \R^n$
    \[ H_P^{(n)}(x) \geq \frac{1}{4}p_4(x)^L > 0.\]
    In particular, $H_P$ is any-dimensional psd.
    \item[(b)] If $P$ has a root in $\N_0^k$, then there are $n \in \N$ and $0 \neq x \in \R^n$ with 
    \[ H_P^{(n)}(x)=-\frac{1}{4}p_4(x)^L < 0.\]
    In particular, $H_P$ is not any-dimensional psd.
\end{enumerate}
Thus, $H_P$ is any-dimensional psd if and only if $P$ has no zero in $\N_0^k$.
\end{proposition}
\begin{proof}

    Consider a basis element $p_{r_1}\cdots p_{r_s}$ of $\Lambda$ occurring in $Q_P$ of degree $e := \sum_{i=1}^s r_i \leq L$. We established $r_j \leq 2m+1$ for $1 \leq j \leq s$ in Remark \ref{rem: S,U,E} and \eqref{eq:degreebound}.
    Under the substitution in \eqref{eq:7} the monomial is mapped to 
    \[ p_4^L \prod_{i=1}^s \frac{p_2^{r_i}p_{2r_i}}{p_4^{r_i}} = p_4^{L-e}p_2^e \prod_{i=1}^s p_{2r_i},\]
    which is a polynomial in $p_2,p_4,\dots,p_{4m+2}$ and homogeneous of degree $4(L-e)+2e+\sum_{i=1}^s 2r_i = 4L$. The constant term $P(0)^2-\frac{1}{4}$ of $Q_P$ is mapped to $\left(P(0)^2-\frac{1}{4}\right)p_4^L$ which is of the same degree. 
    Hence $H_P$ is indeed a homogeneous even symmetric function of degree $4L=4(2(k+2)d+2)$ with rational coefficients.
    Since the constants $C,K_P,L$ and all coefficients of $H_P$ are computed by finite rational arithmetic from $P$, the symmetric function $H_P \in \Lambda_\Q$ is effectively computable.

    We are left with verifying (a) and (b).
    Let $n \in \N$, $0 \neq x \in \R^n$ and let $y$ be as in \eqref{eq:a,b,z}. We obtain from \eqref{eq:7} the identity 
    \[H_P^{(n)}(x)=p_4(x)^L G_P(p_1(y),\dots,p_{2m+1}(y))=p_4(x)^L Q_P(y).\]
    If $P$ has no zero in $\N_0^k$ then $Q_P(w) \geq \frac{1}{4}$ for every $w \in \bigcup_{n \in \N} \R_{\geq 0}^n$ by Lemma \ref{lem:encoding} and so $H_P^{(n)}(x) \geq \frac{1}{4}p_4(x)^L > 0$. Since this holds for all $n$ and $0 \neq x \in \R^n$ this proves (a).
    
    If $P$ has a zero in $\N_0^k$ take $n$ and $z \in \mathcal{A}^n$ from Lemma \ref{lem:encoding} (b) and set $x_i := \sqrt{z_i} \geq 0$. So $p_1(z)=p_2(z)$ and $Q_P(z)=-\frac{1}{4}$. Then, since $z \neq 0$ we have $x\neq 0$ and
    \[ 0 < p_2(x)=\sum_{i=1}^n z_i = p_1(z) = p_2(z)=\sum_{i=1}^n z_i^2 = p_4(x). \]
    Thus $p_2(x)=p_4(x)$ and the normalized vector $\left( \frac{p_2(x)}{p_4(x)} x_1^2,\dots,\frac{p_2(x)}{p_4(x)}x_n^2\right)$ is exactly $z$. Hence Lemma \ref{lem:encoding} implies
    \[ H_P^{(n)}(x)=p_4(x)^L Q_P(z) = -\frac{1}{4}p_4(x)^L < 0\]
    which proves (b).
    
    Finally, note that the claimed equivalence is a direct consequence of the proven statements (a) and (b).
\end{proof}

We can now complete the proof of Theorem \ref{thm:undecidable}.
\begin{proof}[Proof of Theorem \ref{thm:undecidable}]
It follows from Davis, Putnam, Robinson and Matiyasevich's theorem \cite{davis1961decision,matiyasevich1972diophantine} that deciding whether a Diophantine equation $P=0$ for $P \in \Z[Y_1,\dots,Y_k]$ has a solution in $\N_0^k$ is undecidable (see Remark \ref{rem:undecidability}). In \eqref{eq:7} we constructed the effectively computable rational symmetric function $H_P \in \Q[p_2,p_4,\dots,p_{4m+2}]$. All its truncations $H_P^{(n)}$ are polynomials and thus nonnegativity on $\Q^n$ is the same as nonnegativity on $\R^n$ by continuity. It follows from Proposition \ref{prop:undecidable} that $P$ has no zero in $\N_0^k$ if and only if $H_P$ is any-dimensional positive semidefinite. So, the first claim follows.

By \cite{jones1982universal} the decision problem whether a Diophantine equation with nine unknowns has a solution in $\N_0^9$ is undecidable. Thus, Proposition \ref{prop:undecidable} shows that deciding nonnegativity for homogeneous even symmetric functions in $\Q[p_2,p_4,\dots,p_{46}]$ is undecidable, since for $k=9$ we have $m=11$ and so $4m+2=46$.
\end{proof}
\begin{remark}\label{rem:Jin}
Jin's manuscript \cite{Jin2026} was prepared with OpenAI Codex and contains a proof of the undecidability of any-dimensional nonnegativity for symmetric functions in $\Lambda_\Q$ and the author asks for independent human review.
Both proofs encode $N \in \N_0^k$ as multiplicities of coordinates at fixed labels and consider a symmetric approximation. In \cite{Jin2026} no initial restriction to $\bigcup_{n \in \N}\R_{\geq 0}^n$ is required and the obtained symmetric functions are non-homogeneous and not even. Our construction yields homogeneous even symmetric functions and shows that undecidability occurs already for symmetric functions in $\Q[p_2,p_4,\dots,p_{46}]$ via Jones' theorem \cite{jones1982universal}. Moreover, our proof demonstrates that the problem of deciding any-dimensional positive definiteness is also undecidable.   
\end{remark}

\section{A counterexample to Hilbert's 17th problem for symmetric functions}\label{sec:Artin}
Artin's solution \cite{artin1927zerlegung} to Hilbert's 17th problem shows that every nonnegative polynomial is a sum of squares after multiplication by a non-zero sum of squares.
For a nonnegative polynomial $f \in \R[X_1,\dots,X_n]$ we call a sum of squares $0 \neq q \in \R[X_1,\dots,X_n]$ such that $qf$ is a sum of squares a \emph{sos multiplier} of $f$.

In this section we show that an analog of Artin's theorem fails for symmetric functions. In particular, this rules out the existence of an any-dimensional Krivine-Stengle Positivstellensatz.

\begin{theorem}\label{thm:main}
The homogeneous symmetric function
\[F=4p_1^2\left(p_3-4p_1^3\right)^2+9p_2\left(p_6-8p_1^2p_4+16p_1^4p_2\right)-16p_1^8\]
of degree $8$ satisfies the following properties:
\begin{enumerate}
\item $F$ is any-dimensional positive definite, i.e., $F^{(n)}(x)> 0$ for every $x\in\bigcup_{n \in \N}\left(\R^n\setminus\{0\}\right)$.
\item There is no non-zero formal power series $Q \in \R[[X_1,X_2,\dots]]$ of bounded degree such that $Q^{(n)}$ and $Q^{(n)}F^{(n)}$ are sums of squares for every $n \in \N$.
\end{enumerate}
\end{theorem}

In particular, there is no symmetric function $0 \neq Q \in \SSigma$ such that $QF \in \SSigma$. 
Moreover, we show in Proposition \ref{lem:symmetric multiplier} that if $f \in \mathcal{SP}$ has a non-zero sos multiplier in $\R[[X_1,X_2,\dots]]$ of bounded degree then there exists a sos multiplier in $\SSigma$.

We first prove that $F$ is any-dimensional positive definite.
The strategy to prove (2) in Theorem \ref{thm:main} is to show that the degrees of sos multipliers of the truncation $F^{(n)}$ grow at least linearly in $n$.
For this we rely on a theorem which is attributed to Blekherman in \cite{lee2016symmetric}.

\begin{proof}[Proof of Theorem~\ref{thm:main}~(1)]
Consider the function $\sigma: ~ \R \to \Z, ~ t\mapsto\begin{cases}
    0 \hspace{12mm} \text{if }|t|\leq2/3\\
    \operatorname{sgn}(t) \quad \text{else}
\end{cases}$. Then
\begin{align}\label{eq:t sigma}
    \left|t^3-\sigma(t)\right|\leq3t^2\left|t^2-1\right|\qquad(\text{for all }t\in\R).
\end{align}
Indeed, assume without loss of generality $t\geq 0$. The cases $t\in\{0,1\}$ are trivial. If $0<t\leq\frac 2 3$, then $\frac{t^3}{t^2(1-t^2)}=\frac{t}{1-t^2}\leq \frac65 < 3$. If $t>\frac 2 3$ and $t \neq 1$, then
\[
\frac{\left|t^3-\sgn(t)\right|}{t^2\left|t^2-1\right|}=\frac{\left(t^2+t+1\right)\left|t-1\right|}{t^2(t+1)|t-1|}
=\frac1{t^2}+\frac1{t+1}\leq \frac{9}{4}+\frac{3}{5}=\frac{57}{20}<3.
\]
 \\
Now let $x\in\bigcup_{n\in\N}\R^n$ and consider the integer
$k:=\sum_{i} \sigma(x_i)$. 
Write
\[
D:=p_2(x)\left(p_6(x)-2p_4(x)+p_2(x)\right)
=p_2(x)\sum_i x_i^2(x_i^2-1)^2\geq0.
\]
Applying the triangle inequality, the inequality \eqref{eq:t sigma} and Cauchy-Schwarz gives
\begin{equation}\label{eq:round}
\begin{split}
\left|p_3(x)-k\right|&\leq \sum_i \left|x_i^3-\sigma(x_i)\right|
\leq 3\sum_i x_i^2\left|x_i^2-1\right| \leq3\sqrt{\sum_i x_i^2 \sum_i \left(x_i\left|x_i^2-1\right|\right)^2} = 3\sqrt{D}.
\end{split}
\end{equation}
Suppose $x\neq0$. If $p_1(x)=0$, then $F(x)=9p_2(x)p_6(x)>0$.
Otherwise, since $F$ is homogeneous of even degree we may assume $p_1(x)=\frac{1}{2}$. 
Then $F(x)=\left(p_3(x)-\frac{1}{2}\right)^2+9D-\frac{1}{16}$. By \eqref{eq:round} there is an integer $k$ with $|p_3(x)-k|\leq 3\sqrt D$. Hence
\[
\begin{aligned}
\frac{1}{2}&\leq\left|k-\frac{1}{2}\right|= \left|k-p_3(x)+p_3(x)-\frac{1}{2}\right| \leq\Big|p_3(x)-k\Big|+\Big|p_3(x)-\frac{1}{2}\Big|\\
&\leq 3\sqrt D+\left|p_3(x)-\frac{1}{2}\right|
\leq\sqrt{2}\sqrt{\left(p_3(x)-\frac{1}{2}\right)^2+9D}=\sqrt{2}\sqrt{F(x)+\frac{1}{16}},
\end{aligned}
\]
where the last inequality follows from $(a+b)^2 \leq 2(a^2 +b^2)$ for real numbers $a,b$.
Squaring both sides of the previous inequality and subtracting $\frac{1}{16}$ gives $F(x)\geq \frac{1}{16}>0$.
\end{proof}

The proof of Theorem \ref{thm:main} (1) actually shows that $F \geq 16p_1^8$ and $F=9p_2p_6$ if $p_1=0$.

We will need the following theorem from \cite{lee2016symmetric}.\footnote{In \cite{lee2016symmetric} the authors attribute the theorem to Blekherman and refer to his unpublished manuscript \emph{Symmetric sums of squares on the hypercube}.}

\begin{theorem}\cite[Theorem 2.3]{lee2016symmetric}\label{thm:blekherman}
    Let $Q:~ \{0, 1\}^n \to \R_{\geq 0}$ be a symmetric nonnegative real-valued Boolean function and let $q\in\R[T]$ be the unique polynomial of degree at most $n$ such that
\[Q (x_1,\dots, x_n) = q (x_1 + \dots +x_n) \quad \text{for all }x\in\{0,1\}^n.\]
If $Q$ coincides on $\{0,1\}^n$ with a sum of squares of $n$-variate polynomials of degree $d \leq \frac{n}{2}$, then we can write
\[q(T) = \sum_{k=0}^d \sigma_{d-k} \prod_{i=0}^{k-1}(T-i)(n - T-i)\]
with sums of squares $\sigma_k\in \R[T]$ of degree $\deg(\sigma_k) \leq 2k$.
\end{theorem}

\begin{corollary}\label{lem:Boolean-restriction}
Let $n\in\N$ and let
$Q\in\R[X_1,\dots,X_n]$ be a sum of squares with
$\deg Q\le n$. Suppose that $q\in\R[T]$ satisfies
$\deg q\le n$ and
\[
Q(x)=q\left(\sum_{i=1}^n x_i\right)
\qquad\text{for all } x\in\{-1,1\}^n.
\]
Then $q\left(\frac{1}{2}\right)\ge0$.
\end{corollary}
\begin{proof}
Apply Theorem \ref{thm:blekherman} to
\[\widetilde Q: ~ \{0, 1\}^n \to \R_{\geq 0}, ~ x\mapsto Q(2x_1-1,\dots,2x_n-1) \quad \text{and} \quad \widetilde q:=q(2T-n)\in\R[T],\]
i.e. there are sos $\sigma_i\in \R[T]$ with
\[\widetilde q= \sum_{k=0}^{\lfloor n/2\rfloor} \sigma_{{\lfloor n/2\rfloor}-k} \prod_{i=0}^{k-1}(T-i)(n - T-i).\]
Then $q\left(\frac{1}{2}\right)=\widetilde q\left(\frac{n}{2}+\frac{1}{4}\right)\geq 0,$ because
\[\frac{n}{2}+\frac{1}{4}-i>0 \quad \text{and} \quad n-\frac{n}{2}-\frac{1}{4}-i=\frac{n}{2}-\frac{1}{4}-i>0 \quad \text{for } i<\left\lfloor \frac{n}{2}\right\rfloor.\]
\end{proof}

\begin{corollary}\label{cor: bound sos multiplier}
Let $n\in\N$ and $0\ne Q\in\R[X_1,\dots,X_n]$ be a sos multiplier of $F^{(n)}$. Then
\[
\deg Q>\frac{n-8}{2}.
\]
\end{corollary}
\begin{proof}
Write $d:=\deg Q\in2\N_0$, and suppose $n\ge2d+8$.
By applying the Reynolds operator, we can assume that $Q$ is symmetric. 
By the fundamental theorem of symmetric polynomials, we can write
\[Q=\widetilde{Q}\left(p_1^{(n)},\dots,p_d^{(n)}\right) \quad \text{and} \quad F=\widetilde{F}(p_1,\dots,p_6),\]
for some $\widetilde{Q}\in \R[Y_1,\dots,Y_d]$ and $\widetilde{F}\in\R[Y_1,\dots,Y_6]$.
Recall that $d \leq n$ by assumption and so $p_1^{(n)},\dots,p_d^{(n)}$ are indeed algebraically independent.

Set $N=n-d\ge d+8$ and define the univariate polynomials $c_k(T)=\begin{cases}
N,&k\text{ even},\\
T,&k\text{ odd}.
\end{cases}$

For $z\in \R^d$, define the univariate polynomials
\[
q(T,z)=\widetilde{Q}\bigl(c_1(T)+p_1^{(d)}(z),\dots,c_d(T)+p_d^{(d)}(z)\bigr) ~ \text{and} ~
f(T,z)=\widetilde{F}\bigl(c_1(T)+p_1^{(d)}(z),\dots,c_6(T)+p_6^{(d)}(z)\bigr).
\]
For $x\in\{-1,1\}^N$ and $z\in \R^d$ we have $p_k^{(n)}(x,z)=c_k\left(\sum_{i=1}^N x_i\right)+p_k^{(d)}(z)$ and therefore
\[
Q(x,z)=q\left(\sum_{i=1}^N x_i,z\right) \quad \text{and} \quad (QF^{(n)})(x,z)=(qf)\left(\sum_{i=1}^N x_i,z\right).
\]
Note that, for fixed $z \in \R^d$, the left-hand sides are sos polynomials in $x$.
Moreover, they and the corresponding polynomials in the variable $T$
have degrees at most $d+8\le N$. By Corollary~\ref{lem:Boolean-restriction},
\begin{equation}\label{eq:Boolean-nonnegative}
q\left(\frac{1}{2},z\right)\ge0 \quad \text{and} \quad q\left(\frac{1}{2},z\right)f\left(\frac{1}{2},z\right)\ge0.
\end{equation}
At $z=0$, we have $c_k\left(\frac{1}{2}\right)+p_k^{(d)}(0)=c_k\left(\frac{1}{2}\right)=\begin{cases}
    N  \quad \text{for }k\in\{2,4,6\}\\
    \frac12 \quad \hspace{1.3mm} \text{for }k\in\{1,3\}
\end{cases}$ and therefore 
\[
f\left(\frac{1}{2},0\right)=4\left(\frac{1}{2}\right)^2\Bigg(\underbrace{\left(\frac{1}{2}\right)-4\left(\frac{1}{2}\right)^3}_{=0}\Bigg)^2+9N\Bigg(\underbrace{N-8\left(\frac{1}{2}\right)^2N+16\left(\frac{1}{2}\right)^4N}_{=N-2N+N=0}\Bigg)-16\left(\frac{1}{2}\right)^8=-\frac{1}{16}<0 .
\]
Consequently $f\left(\frac{1}{2},0\right)<0$ on a nonempty open neighborhood of $0$ and therefore $q\left(\frac{1}{2},z\right)=0$ on this neighborhood by Equation~\eqref{eq:Boolean-nonnegative}. Hence, $q\left(\frac{1}{2},\underline{Z}\right)\in\R[Z_1,\dots,Z_d]$ is the zero polynomial.
This contradicts $\widetilde{Q}\ne0$: the polynomials
$p_1^{(d)},\dots,p_d^{(d)}$ are algebraically
independent, as are their translates $c_k(1/2)+p_k^{(d)}$.
Thus $\deg Q>\frac{n-8}{2}$.
\end{proof}
We can now finish the proof of Theorem \ref{thm:main}. 
\begin{proof}[Proof of Theorem~\ref{thm:main} (2)]
Suppose that a non-zero formal power series $Q$ of
bounded degree satisfies $Q^{(n)}$ and $Q^{(n)}F^{(n)}$
are sos for all $n$. Choose an integer $D$ with $\deg Q^{(n)}\le D$
for every $n$. Because $Q\ne0$, one of its monomials with non-zero coefficient involves only finitely many variables, and hence
$Q^{(n)}\ne0$ for all sufficiently large $n$.
Corollary~\ref{cor: bound sos multiplier} then gives
\[
D\ge\deg Q^{(n)}>\frac{n-8}{2}
\]
for all sufficiently large $n$, a contradiction.
\end{proof}

\begin{remark}
The proof of Theorem \ref{thm:main} also shows why Artin's solution to Hilbert's 17th problem does not extend to symmetric functions: Nonnegativity on $\bigcup_{n\in\N}\R^n$ does not imply nonnegativity on the real spectrum
\[\sper_{\SSigma}\left(\Lambda\right)= \left\{ P ~ \middle| ~ P \text{ ordering of $\Lambda$ with }\SSigma\subseteq P \right\}.\]
To see this, consider the order $\R(T)_{\geq 0}$ on $\R(T)$, where $T\geq c$ for all $c\in \R$. Then the homomorphism
\[\varphi: ~ \Lambda \to \R(T), ~ p_k\mapsto\begin{cases}
    T \quad \text{for $k$ even}, \\
    \frac{1}{2} \quad \text{for $k$ odd}
\end{cases}\]
gives rise to an ordering $P:=\varphi^{-1}\left(\R(T)_{\geq 0}\right)$. Corollary \ref{lem:Boolean-restriction} shows that $\SSigma \subseteq P$, so $P\in \sper_{\SSigma}\left(\Lambda\right)$ and the proof of Corollary \ref{cor: bound sos multiplier} shows that $\varphi(F)=-\frac{1}{16}\notin\R(T)_{\geq 0}$, so $F\notin P$.
\end{remark}

\begin{corollary}\label{cor:no artin analog}
    There exist symmetric functions that are any-dimensional positive definite which are not sums of squares in the quotient field of $\Lambda$. Moreover, no any-dimensional Krivine-Stengle Positivstellensatz exists.
\end{corollary}
\begin{proof}
    The first statement follows immediately from Theorem \ref{thm:main} since every symmetric function is a formal power series of bounded degree.
    Suppose that a Krivine-Stengle Positivstellensatz for any-dimensional nonnegative symmetric functions exists. Then for every $f \in \mathcal{SP}$ there are some $0 \neq q,g \in \SSigma$ and $k \in \N_0$ with $qf=f^{2k}+g \in \SSigma$. However, $F$ from Theorem \ref{thm:main} serves as a counterexample.
\end{proof}

\begin{remark}
In \cite[Theorem~5.4]{debus2026any} the authors proved that if a homogeneous symmetric function $G \in \Lambda$ written as $G=\sum_{i=0}^{2d} G_i p_1^{i}$ ($G_i\in \R[p_2,p_3,p_4,\dots]$, $G_{2d}\neq 0$) satisfies 
\[\inf_{p_2^{(n)}=1} G^{(n)} > \varepsilon \quad \text{and} \quad \inf_{p_2^{(n)}=1} G_{2d}^{(n)} > \varepsilon \quad \text{for all $n \in \N$ and some $\varepsilon > 0$}\]
then there exists $r \in \N$ with $p_2^r G \in \S\Sigma$. The second condition is fulfilled for $F$, but the first condition fails: Although all truncations of the octic $F$ from Theorem \ref{thm:main} are positive definite their minima on the unit sphere are not bounded away from zero uniformly in $n$.
Note that at $x = (-1,\dots,-1,1,\dots,1) \in \R^{2n}$ with $n$ many $-1$'s and $1$'s we have
$F^{(2n)}(x) =  36 n^2 $. Since $p_2(x)=2n$ and $F$ is homogeneous of degree $8$, $F$ evaluates at the corresponding point on the sphere $\frac{1}{\sqrt{2n}}x$ to $\frac{36n^2}{(2n)^{4}}=\frac{9}{4}n^{-2}.$
In particular, Theorem \ref{thm:main} shows that the assumption in \cite[Theorem~5.4]{debus2026any} cannot be weakened to any-dimensional positive definiteness.    
\end{remark}

\begin{remark}
Blekherman, Gouveia and Pfeiffer construct for each $n$ a globally positive symmetric polynomial $f_n \in \R[X_1,\dots,X_n]$ of degree $4$ which only has sos multipliers of degree at least $2\left\lfloor \frac{n}{2} \right\rfloor -2$ (see \cite[Theorem~4.1]{BlekhermanDegree}). However, the coefficients of $f_n$ depend on $n$ and the $f_n$ are not the truncations of the same symmetric function. Their proof uses sums of squares on the Boolean hypercube and the representation theory of the symmetric group. Corollary \ref{cor: bound sos multiplier} is a similar statement with the weaker lower bound $n/2-4$, but there the polynomials $F^{(n)}$ are the truncations of one fixed symmetric function $F$ (of degree $8$ rather than $4$). 
\end{remark}

For a partition $\lambda = (\lambda_1,\dots,\lambda_\ell) \vdash d$ we write $p_\lambda=p_{\lambda_1}\dots p_{\lambda_\ell}$.
Finally, we show that if a symmetric function in $\mathcal{SP}$ admits a non-zero sos multiplier in $\R[[X_1,X_2,\dots]]$ of bounded degree then it admits a sos multiplier in $\SSigma$.

\begin{proposition}\label{lem:symmetric multiplier}
    Let $f \in \mathcal{SP}$ be any-dimensional psd. Suppose that there exist $d \in \N$ and $0 \neq q \in \R[[X_1,X_2,\dots]]$ with $\deg q^{(n)} \leq 2d$ and both polynomials $q^{(n)}$ and $q^{(n)}f^{(n)}$ are sos for every $n \in \N$. Then there exists an any-dimensional sos $ h \in \SSigma \setminus \{0\}$ with $hf \in \SSigma$.
\end{proposition}
\begin{proof}
For each $n \in \N$ we can symmetrize $q^{(n)}$. Thus, there exists a sequence 
\[h_n := \frac{1}{n!} \sum_{\sigma \in S_n} q^{(n)}\left(X_{\sigma(1)},\dots,X_{\sigma(n)}\right) \]
of symmetric sos polynomials with $h_nf^{(n)}$ sos. 
Note that $h_n$ is non-zero if $q^{(n)} \neq 0$, since $q^{(n)}$ is a sum of squares. Since $q \neq 0$ we have $q^{(n)} \neq 0$ for all sufficiently large $n$ and thus also $h_n \neq 0$. 
Write $h_n=\sum_{\substack{\lambda \vdash r  \\ r \leq 2d}} c_{n,\lambda} p_\lambda^{(n)}$ in the power sum basis for $n \geq 2d$. Rescaling $h_n$ by a positive scalar we can assume that $\sum_{\substack{\lambda \vdash r  \\ r \leq 2d}} c_{n,\lambda}^2=1$ for sufficiently large $n$. 
The associated sequence of symmetric functions $\Tilde{h}_n := \sum_{\substack{\lambda \vdash r  \\ r \leq 2d}} c_{n,\lambda} p_\lambda \in \Lambda$ must have a convergent subsequence $\Tilde{h}_{n_k}$ with limit $0 \neq h\in \Lambda$. The convergence is considered coefficient-wise in the power sum basis. 
For each fixed $m \in \N$ and $n_k$ sufficiently large the symmetric polynomial $\Tilde{h}_{n_k}^{(m)}$ is sos, since $h_{n_k} \in \R[X_1,\dots,X_{n_k}]$ is sos and thus also all its truncations. The sequence $\Tilde{h}_{n_k}^{(m)}$ converges to $h^{(m)}$. Since the cone of sums of squares of degree $\leq 2e$ in $\R[X_1,\dots,X_m]_{\leq 2e}$ is closed \cite[Proposition 2.4.6]{Scheiderer} for every $e \in \N$, the polynomials $h^{(m)}$ and $h^{(m)}f^{(m)}$ are sos. 
Because $m$ was arbitrary the claim holds for all $m \in \N$ and thus $h \in \Lambda$ and $hf \in \Lambda$ are both any-dimensional sos.    
\end{proof}

\section{Dimension-free trace inequalities and entanglement criteria}\label{sec:quantum}
In this section we follow the notation and exposition in \cite{mallick2026optimal}, where, however, the degree is denoted by $n$ and the matrix size by $d$. To align with Sections \ref{sec:undecidability} and \ref{sec:Artin} we use $n$ to denote the matrix size and $d$ for the degree.

Quantum entanglement is the defining characteristic which distinguishes quantum correlations from their classical counterparts. In \emph{entanglement theory} the central computational challenge is known as the \emph{separability problem}: determining whether a given bipartite density matrix $\rho$ acting on the joint Hilbert space $\C^{d_A} \otimes \C^{d_B}$ is separable. 
Physically, a state is separable if it can be prepared using only local operations and classical communication. Mathematically, a state $\rho$ is \emph{separable} if it can be decomposed into a conic combination of local product states, i.e.,
\[ \rho = \sum_i a_i \cdot  \sigma_A^i \otimes \sigma_B^i\]
where $a_i \geq 0$ and the local states $\sigma_A^i$ and $\sigma_B^i$ are positive semidefinite. Since this problem is NP-hard with respect to the local dimensions $d_A$ and $d_B$ \cite{gurvits2003classical}, the underlying computations become intractable for large quantum systems. Consequently, the practical framework is built upon separability criteria rather than exact solutions. These criteria are computationally efficient necessary conditions for a state to be separable. 

Randomized measurement protocols give experimental access to polynomial invariants of a quantum state $\rho$ such as the moments $\tr{(\rho^k)}$ and the moments $\tr((\rho^{\Gamma_B})^k)$ of its partial transpose $\rho^{\Gamma_B}$ without tomography of the state itself \cite{elben2020mixed, yu2021optimal}.
Entanglement criteria relying on such invariants are polynomial inequalities which hold for all separable states. To be applicable without knowledge of the local dimensions these inequalities must hold in any dimension.

Mallick, Gulati and Nechita \cite{mallick2026optimal} study the convex cones of such dimension-free inequalities. Their one-matrix versions are the convex cones of (even) homogeneous symmetric functions that are any-dimensional psd and were previously studied in \cite{acevedo2025symmetric, acevedo2025wonderful,levin2025any}. 

We recall the setup from \cite{mallick2026optimal} and mention the consequences of Theorems \ref{thm:undecidable} and \ref{thm:main} for dimension-independent entanglement criteria. \\

For a positive integer $n$ let $\mathcal{M}_n$ denote the space of complex $n \times n$ matrices, $\mathcal{M}_n^{\operatorname{sa}}$ the space of self-adjoint or Hermitian matrices in $\mathcal{M}_n$ and $\mathcal{M}_n^+$ the cone of positive semidefinite matrices in $\mathcal{M}_n$. 
By the spectral theorem all eigenvalues of $X \in \mathcal{M}_n^{\operatorname{sa}}$ are real. Also, all eigenvalues of matrices in $\mathcal{M}_n^+$ are nonnegative.  
Let $X \in \mathcal{M}_n^{\operatorname{sa}}$ and $x \in \R^n$ be a vector of eigenvalues of $X$ (which is unique up to permutation). For $r \in \N$ we define 
\[p_r(X):=  \tr (X^r)=p_r(x)\]
which is well-defined, since the power sums are symmetric polynomials and accordingly $f(X):=f^{(n)}(x)$ for $f \in \Lambda$. 
The linear span of all $p_\lambda$ with $\lambda \vdash d$ over $\R$ is the space $\Lambda_d \subset \Lambda$ of homogeneous symmetric functions of degree $d$.

In \cite{mallick2026optimal} the following \emph{dimension-free} cones are studied:
\begin{align*}
    &\mathsf{Ineq}^{\mathrm{all}}_{d} ~:= \left\{f \in \Lambda_{d} ~\middle|~ f(X) \geq 0 \text{ for all } n \geq 1 \text{ and all } X \in \mathcal{M}_n^{\mathrm{sa}}\right\} = \mathcal{SP} \cap \Lambda_d, \\
    &\mathsf{Ineq}^{\mathrm{psd}}_{d} := \left\{f \in \Lambda_{d} ~\middle|~ f(X) \geq 0 \text{ for all } n \geq 1 \text{ and all } X \in \mathcal{M}_n^{\mathrm{+}}\right\} ~\cong \mathcal{SP} \cap \Lambda_{2d}^{\varepsilon},
\end{align*}
where $\Lambda_{2d}^{\varepsilon}$ denotes the even symmetric functions in $\Lambda_{2d}$ and the isomorphism is $f \mapsto f\left(X_1^2,X_2^2,\dots\right)$, i.e., the substitution $p_r \mapsto p_{2r}$. 
The dimension-free cones come together with the moment cones $\mathsf{Mom}^{\mathrm{all}}_d$ and $\mathsf{Mom}^{\mathrm{psd}}_d$ which are the closed convex conical hulls of all vectors $(p_\lambda(X))_{\lambda \vdash d}$ for $X \in \bigcup_{n \in \N}\mathcal{M}_n^{\mathrm{sa}}$ and $X \in \bigcup_{n \in \N}\mathcal{M}_n^{+}$ respectively. So, by construction the moment cones are the dual cones of the dimension-free cones (\cite[Proposition~2.3]{mallick2026optimal}).
Non-semialgebraicity and thus non-polyhedrality of the cones $\mathsf{Mom}^{\mathrm{all}}_d$ and $\mathsf{Ineq}^{\mathrm{all}}_{d}$ for all even $d \geq 4$, and of the cones $\mathsf{Mom}^{\mathrm{psd}}_d$ and $\mathsf{Ineq}^{\mathrm{psd}}_{d}$ for all $d \geq 3$ was proven in \cite[Theorem~3.12]{acevedo2025symmetric}. The extreme rays of $\mathsf{Ineq}^{\mathrm{psd}}_{3}$ and $\mathsf{Mom}^{\mathrm{psd}}_{3}$ were characterized in \cite{Choi} (see also \cite[Proposition~3.5]{mallick2026optimal}). Moreover, the slice $p_1=1$ of $\mathsf{Mom}^{\mathrm{psd}}_d$ is the closed convex hull of the image of the Vandermonde cell $\Pi_d$ from \cite{acevedo2025wonderful} under $(p_2,\dots,p_d) \mapsto (p_\lambda)_{\lambda \vdash d}$ where $p_1$ is set to $1$. 
The Vandermonde cell $\Pi_d$ itself is not convex for $d\geq 3$ and has the combinatorial structure of an infinite cyclic polytope \cite[Theorem~4.4]{acevedo2025wonderful}.
For $d=3$ the slice of $\mathsf{Mom}^{\mathrm{psd}}_3$ is $\operatorname{conv} \Pi_3$. The boundary of $\Pi_3$ was parametrized in \cite[Corollary 2.22]{acevedo2025wonderful} (see also Figure 2 and Corollary 3.7 in \cite{mallick2026optimal}).
However, non-semialgebraicity does not rule out the existence of effective certificates for membership in $\mathsf{Ineq}^{\mathrm{all}}_{d}$ and $\mathsf{Ineq}^{\mathrm{psd}}_{d}$.

For $w = \sum_{\sigma,\tau \in S_d} c_{\sigma,\tau} (\sigma,\tau)$ in the group algebra $\R[S_d\times S_d]$ and a matrix $\rho$ acting on $\C^{d_A} \otimes \C^{d_B}$ the \emph{trace invariant} of $\rho$ with respect to $w$ is
\[ \tr_w(\rho):= \operatorname{Re} \sum_{\sigma,\tau \in S_d} c_{\sigma,\tau} \tr \left(\left(P_\sigma \otimes P_\tau\right) \rho^{\otimes d}\right),\]
where $P_\sigma$ and $P_\tau$ permute the $d$ tensor copies of $\C^{d_A}$ and $\C^{d_B}$ according to $\sigma$ and $\tau$ respectively.

The \emph{dimension-free separability-valid trace-invariant cone} of \cite{mallick2026optimal} is
\[ \mathsf{Ineq}^{\mathrm{sep}}_{d} ~:= \left\{w \in \R[S_d \times S_d] ~\middle| ~ \tr_w(\rho) \geq 0 \text{ for every $d_A, d_B \geq 1$ and every $\rho \in \mathrm{Sep}\left(\C^{d_A} : \C^{d_B}\right)$}\right\}, \]
where for local dimensions $d_A,d_B$ the set
\[  \mathrm{Sep}\left(\C^{d_A} : \C^{d_B}\right) := \operatorname{cone}\left\{X \otimes Y ~\middle|~ X \in \mathcal{M}_{d_A}^+, Y \in \mathcal{M}_{d_B}^+\right\}\]
is the \emph{unnormalized separable cone}. 
If $\alpha_\lambda \in S_d$ denotes a permutation of cycle type $\lambda$ then 
\[\tr_{(\alpha_\lambda,\alpha_\lambda)}(\rho)=p_\lambda(\rho)\quad \text{and}\quad \tr_{(\alpha_\lambda,\alpha_\lambda^{-1})}(\rho)=p_\lambda\left(\rho^{\Gamma_B}\right).\]
In \cite[Proposition 4.3]{mallick2026optimal} the cones $\mathsf{Ineq}_d^{\mathrm{all}} \subset \mathsf{Ineq}_d^{\mathrm{psd}}$ are embedded into $\mathsf{Ineq}_d^{\mathrm{sep}}$.
Since $\rho$ and $\rho^{\Gamma_B}$ are positive semidefinite whenever $\rho$ is separable the  linear maps
\[ \iota_+(f) := \sum_{\lambda \vdash d}a_\lambda (\alpha_\lambda,\alpha_\lambda) \quad \text{and} \quad \iota_\Gamma(f) := \sum_{\lambda \vdash d}a_\lambda \left(\alpha_\lambda,\alpha_\lambda^{-1}\right) \text{ for } f=\sum_{\lambda \vdash d}a_\lambda p_\lambda \in \Lambda_d\]
embed the one-matrix cone $\mathsf{Ineq}^{\mathrm{psd}}_{d}\supset \mathsf{Ineq}^{\mathrm{all}}_{d}$ into $\mathsf{Ineq}^{\mathrm{sep}}_{d}$. 
Conversely, for $X \in \mathcal{M}_n^+$ the matrix $\rho := X \otimes 1 \in \mathrm{Sep}\left(\C^n : \C^1\right)$ satisfies $\rho = \rho^{\Gamma_B}=X$ and therefore
\begin{align}\label{eq:equivalences}
    f \in \mathsf{Ineq}^{\mathrm{psd}}_{d} \iff \iota_+(f) \in \mathsf{Ineq}^{\mathrm{sep}}_{d} \iff \iota_\Gamma(f) \in \mathsf{Ineq}^{\mathrm{sep}}_{d}.
\end{align}
The image of $\iota_\Gamma$ consists of the dimension-free criteria in terms of the moments of the partial transpose which form the framework of \cite{elben2020mixed,yu2021optimal}.

\begin{remark}
The trace in the definition of $\mathsf{Ineq}^{\mathrm{all}}_{d}, \mathsf{Ineq}^{\mathrm{psd}}_{d}$ and $\mathsf{Ineq}^{\mathrm{sep}}_{d}$ is not normalized. A density matrix satisfies $\tr \rho = 1$ and not $\frac{1}{n} \tr \rho = 1$. 
The setting of trace polynomials with normalized trace $\frac{1}{n}\tr$ in \cite{Klep26,Klep21} corresponds to the power mean inequalities mentioned in the introduction which were studied in \cite{AcevedoBlekherman,BlekhermanRiener}. In this setting any-dimensional nonnegativity is decidable and sos Positivstellensätze exist \cite{Klep26,Klep21}. Theorem \ref{thm:quantum1} shows that the dimension-free cones from \cite{mallick2026optimal} behave fundamentally differently.
\end{remark}

We can reformulate Theorem \ref{thm:undecidable} as assertions on quantum entanglement showing that membership in the one-matrix cones $\mathsf{Ineq}^{\mathrm{all}}_{d}$ and $\mathsf{Ineq}^{\mathrm{psd}}_{d}$ cannot be decided algorithmically.

\begin{theorem}\label{thm:quantum1}
The following decision problems are undecidable:
\begin{enumerate}
    \item[(a)] Instance: A homogeneous $f \in \Lambda_\Q$. Question: Is $f \in \mathsf{Ineq}^{\mathrm{all}}_{\deg f}$? \\
    The problem remains undecidable when the instances are restricted to $f \in \Q[p_2,p_4,\dots,p_{46}]$.
    \item[(b)] Instance: A homogeneous $f \in \Lambda_\Q$. Question: Is $f \in \mathsf{Ineq}^{\mathrm{psd}}_{\deg f}$? \\
    The problem remains undecidable when the instances are restricted to $f \in \Q[p_1,p_2,\dots,p_{23}]$.
    \item[(c)] Instance: A $d \in \N$ and $w \in \Q[S_d \times S_d]$. Question: Is $w \in \mathsf{Ineq}_d^{\mathrm{sep}}$? \\
    The problem remains undecidable when the instances are restricted to the image of $\iota_+$ or to the image of $\iota_\Gamma$.
\end{enumerate}
In particular, there is no algorithm deciding whether a polynomial inequality with rational coefficients in $\tr(\rho),\dots,\tr(\rho^{23})$ holds for all density matrices of all dimensions. Moreover, there is no algorithm deciding whether a polynomial inequality with rational coefficients in the moments $\tr\left(\rho^{\Gamma_B}\right),\dots,\tr\left(\left(\rho^{\Gamma_B}\right)^{23}\right)$ of the partial transpose holds for all separable states in all local dimensions.
\end{theorem}
\begin{proof}
Assertion (a) is a reformulation of Theorem \ref{thm:undecidable}. 

For (b) note that $f \in \mathsf{Ineq}^{\mathrm{psd}}_d$ if and only if $f\left(X_1^2,X_2^2,\dots\right) \in \mathsf{Ineq}^{\mathrm{all}}_{2d}$ and then the statement follows from (a).

For (c) we use \eqref{eq:equivalences}. Since the linear maps $\iota_+$ and $\iota_\Gamma$ from \cite{mallick2026optimal} are effectively computable the statement follows from (b).

To see the last claims we observe that a homogeneous symmetric function $f \in \Lambda_d$ is nonnegative on $\bigcup_{n \in \N} \mathcal{M}_n^+$ if and only if $f$ is nonnegative on all density matrices, since $f$ is homogeneous and every non-zero matrix $X \in \mathcal{M}_n^+$ is a positive multiple of a density matrix.
A polynomial inequality in $\tr (\rho),\dots,\tr(\rho^{23})$ is a polynomial in $p_1,\dots,p_{23}$ evaluated at the eigenvalues of $\rho$. So the instances of (b) with $f \in \Q[p_1,\dots,p_{23}]$ form a subset of the instances. The same argument can be applied to $\rho^{\Gamma_B}$ instead of $\rho$. 
\end{proof}

\begin{remark}
    It is worth comparing the NP-hardness of the separability problem \cite{gurvits2003classical} which applies for a given state in fixed local dimensions with Theorem \ref{thm:quantum1}. The problem of verifying the validity of polynomial trace-invariant inequalities for all separable states in fixed local dimensions $d_A,d_B$ is decidable by the Tarski-Seidenberg theorem \cite{Seidenberg,Tarski} since $\mathrm{Sep}\left(\C^{d_A} : \C^{d_B}\right)$ is a semialgebraic set. The universality over all dimensions makes the validity problem of Theorem \ref{thm:quantum1} undecidable. 
\end{remark}

\begin{proposition}\label{prop:quantum1}
    The symmetric function $F$ of Theorem \ref{thm:main} lies in $\mathsf{Ineq}^{\mathrm{all}}_8 \subset \mathsf{Ineq}^{\mathrm{psd}}_8$ and is positive on every non-zero Hermitian matrix. 
    There is no non-zero $Q \in \R[[X_1,X_2,\dots]]$ of bounded degree such that $Q^{(n)}$ and $(QF)^{(n)}$ are sos for every $n \in \N$, i.e., $F$ admits no sums of squares certificate that is uniform in the dimension, even if sum of squares denominators are allowed.
\end{proposition}
\begin{proof}
    This follows directly from Theorem \ref{thm:main}.
\end{proof}

\subsection*{Acknowledgments}
The authors thank Jose Acevedo, Greg Blekherman and Cordian Riener for earlier discussions on the undecidability question.

\subsection*{AI Declaration}
This work relies on significant input from large language models. Initial proofs of the two main results (Theorem \ref{thm:undecidable} and Theorem \ref{thm:main}) were found by GPT 6.0 Astra (OpenAI) in September 2026. The authors then verified, simplified and strengthened the proofs, with additional assistance from Claude Fable 5.1 (Anthropic). The final presentation was written by the authors and AI was used to
improve grammar, spelling and readability. The authors take full responsibility for the content.\\  
Finally, a Lean formalization of our proof was almost entirely generated by coding agents. The formalization is available on GitHub at
\begin{center}
    \url{https://github.com/RobinSchabert/positive-sym-func-h17-undecidability-lean}.
\end{center}
\bibliographystyle{abbrv}
\bibliography{references}

\end{document}